\documentclass[12pt,oneside]{article}
\usepackage{amsmath,amsthm,amsfonts,amssymb,MnSymbol,mathrsfs,latexsym}
\usepackage[english]{babel}
\usepackage{bm}

\numberwithin{equation}{section}

 \newcommand{\const}{\rm const}

\theoremstyle{plain}
\newtheorem{theorem}{Theorem}[section]
\theoremstyle{theorem}

\newtheorem{proposition}[theorem]{Proposition}

\newtheorem{remark}{Remark}[section]
\newtheorem{example}{Example}[section]

\renewenvironment{proof}{{\bf{Proof.}}}{\hfill $\Box$ \\}

\DeclareMathOperator*{\esssup}{ess\,sup}

\title{\large \textbf{Quasi Grand Lebesgue Spaces}}

\footnotesize\date{}

\author{\normalsize Maria Rosaria FORMICA ${}^{1}$,   \normalsize Eugeny OSTROVSKY
${}^2$ and \normalsize Leonid SIROTA ${}^3$}

\begin{document}

\maketitle

\begin{center}
{\footnotesize ${}^{1}$ Universit\`{a} degli Studi di Napoli \lq\lq Parthenope\rq\rq, via Generale Parisi 13,\\
Palazzo Pacanowsky, 80132,
Napoli, Italy.} \\

\vspace{2mm}

{\footnotesize e-mail: mara.formica@uniparthenope.it} \\

\vspace{4mm}

{\footnotesize ${}^{2,\, 3}$  Bar-Ilan University, Department of Mathematics and Statistics, \\
52900, Ramat Gan, Israel.} \\

\vspace{2mm}

{\footnotesize e-mail: eugostrovsky@list.ru}\\

\vspace{2mm}

{\footnotesize e-mail: sirota3@bezeqint.net} \\

\end{center}

 \vspace{3mm}

\begin{abstract}
 We introduce a new class of quasi-Banach spaces as an extension of the classical Grand Lebesgue Spaces for \lq\lq small\rq\rq values
 of the parameter, and we investigate some its properties, in particular, completeness, fundamental function, operators estimates, Boyd indices,
 contraction principle, tail behavior, dual space, generalized triangle and quadrilateral constants and inequalities.
 \end{abstract}

 \vspace{3mm}


\vspace{2mm}

\noindent {\it \footnotesize 2010 Mathematics Subject
Classification}:
{\footnotesize  46E30,   % Spaces of measurable functions
46A16 % Not locally convex spaces (metrizable topological linear spaces, locally bounded spaces, quasi-Banach spaces, etc.)
% 46B10    % Duality and reflexivity in normed linear and Banach spaces
% 60B05    % Probability measures on topological spaces
% 26D15    % Inequalities for sums, series and integrals
}

 \vspace{4mm}

 \section{Introduction}

\vspace{2mm}

Let $X$ be a vector (linear) space. A function $||\cdot||:
X\to[0,\infty)$ is said a \emph{quasi-norm} if

\begin{equation} \label{axiom 1}
||x|| \geq 0, \ \ \forall x \in X  ; \  \ \ ||x|| = 0 \
\Leftrightarrow \ x = 0;
\end{equation}

\begin{equation} \label{axiom 3}
||\alpha \, x|| = |\alpha| \cdot ||x||, \ \ \ \forall x \in X,  \
\alpha \in \mathbb R.
\end{equation}

\begin{equation} \label{axiom 2}
\exists \, C \in [1,\infty) \ : \  \ ||x+y||  \le C(||x|| + ||y||),
\ \ \ \ \forall x,y \in X
\end{equation}

The space $X$, equipped with the quasi-norm $||\cdot||=||\cdot||_X$,
is called \emph{quasi-normed space}. Of course, if $C=1$ in
\eqref{axiom 2}, then $||\cdot||$ is a norm, i.e. $X$ is a normed
space.

\vspace{2mm}

This definition was introduced by Tosio Aoki in \cite{Aoki}; see
also \cite{Burenkov1,Hudzikandt,Kalton,Maligranda 1,Maligranda
2,Mitani,Rolewicz,Wu}.

 We highlight the paper of Lech Maligranda \cite{Maligranda 1}, which will very useful for our purposes. A modern survey of the theory of these
 spaces, containing  new results and applications, may be found in the recent
article \cite{Nekvinda} (see also reference therein).

\vspace{4mm}

 The smallest possible constant  $ \ C = C(X) \ $ in (\ref{axiom 2}) is called the \emph{quasi-triangle constant} of the space
$ \  X = (X,||\cdot||). \ $ \par

 If, in addition, we have for some (constant) value $p\in(0,1]$,

\begin{equation} \label{p norm}
||x + y||^p \le ||x||^p + ||y||^p,
\end{equation}
then the functional $ \ x \to ||x|| \ $ is called a {\it p-norm} and
the value $p$ is  named as {\it power parameter}.

\vspace{3mm}

It is known, thanks to Aoki-Rolewicz theorem (see
\cite{Aoki,Rolewicz}), that if $ \ (X,||\cdot||) \ $ is
 a complete quasi-normed space, then there exists a constant $p\in(0,1)$  and a $p$\,-norm $ ||x||' $ in $X$ equivalent
 to the source one $ ||\cdot||$. \par
 More precisely, if  $\,  C\, $ is the quasi-triangle constant, then  $ p$ may be determined by the relation $C = 2^{1/p - 1} $ and
$$
||x||' \le ||x|| \le 2^{1/p} ||x||', \ \ \ \forall \, x \in X,
$$
 (see also \cite{Gustavsson}).

 % (see also The Extension of the Finite-Dimensional Version
% of Krivineís Theorem to Quasi-Normed Spaces ALEXANDER E. LITVAK, Lemma 1)

 \vspace{2mm}

Let now and in the sequel $ \ (Q, B,\mu) \ $ be a non-trivial
measure space with sigma-finite measure $\mu$.

 \vspace{2mm}

 An important class of quasi-Banach spaces which are not Banach
spaces is the class of $L_p$ spaces, for $ \ 0 < p < 1 \ $, with the
usual quasi-norm

\begin{equation} \label{Lp example}
||f||_p :=  \left[ \int_Q   \ |f(x)|^p \ \mu(dx) \ \right]^{1/p}.
\end{equation}

In this case
\begin{equation} \label{Lp integral estim}
||f + g||_p \le 2^{1/p -1} (||f||_p + ||g||_p),
\end{equation}
i.e., the quasi-triangle constant of $ \, L_p \, $ is $\,2^{1/p
-1}\,$ (see \cite{Maligranda 1}). \par

 \vspace{1mm}

 Note, in addition, that in this case
\begin{equation} \label{p norm}
 ||f+g||_p^p \le ||f||_p^p + ||g||_p^p.
 \end{equation}

 Both the estimates (\ref{Lp integral estim}) and (\ref{p norm}) follows immediately  from the double inequality
\begin{equation} \label{ineq A B}
(A + B)^p \le    A^p + B^p \le 2^{1 - p} (A + B)^p, \ \ \ A, B \ge
0, \ \ \ 0<p<1.
\end{equation}

\vspace{3mm}

 We introduce in this paper a new class of quasi-Banach spaces, as
 an extension of the known Grand Lebesgue Spaces, but for \lq\lq small\rq\rq values of the degree
parameter (see Section \ref{Section Quasi Grand}), and we
investigate some its properties.

 \vspace{4mm}

 \section{Ordinary Grand Lebesgue Spaces}\label{Section_OrdinaryGLS}

\vspace{2mm}

 We recall here briefly the definition of (ordinary) Grand Lebesgue Space (GLS).

\vspace{2mm}

Let $ (Q = \{ x \}, B,\mu)$  be a measure space with non-trivial
sigma-finite measure $\mu$ and let
 $\psi(p) = \psi_{\alpha,\beta}(p), \ p \in (\alpha,\beta)$,
where $ 1 \le \alpha < \beta \le \infty $, be a continuous, strictly
positive, numerical valued function in the {\it open} interval such
that
$$
\inf_{p \in (\alpha, \beta)} \psi(p) > 0.
$$

The Grand Lebesgue Space (GLS) \, $ G\psi = G\psi_{\alpha,\beta}$ is
a rearrangement invariant Banach function space in the classical
sense (\cite[Chapters 1,2]{Bennet}) and consists of all the
measurable functions $ f: Q \to \mathbb R$ having finite norm
\begin{equation} \label{Grand Lebesgue}
||f||_{G\psi} = ||f||_{G\psi_{\alpha,\beta}}  \stackrel{def}{=}
\sup_{p \in (\alpha,\beta)} \left\{ \ \frac{||f||_{p}}{\psi(p)} \
\right\}.
\end{equation}

\begin{remark}
{\rm Define $\psi(p)=(\beta-p)^{-\frac{\theta}{p}}, \ p \in
(1,\beta), \ 1 < \beta < \infty, \ \theta\geq0$. Let $Q\subset
\mathbb R^n$, \, $n\geq1$, be a measurable set with finite Lebesgue
measure.

\noindent Replace $p$ with $\beta-\varepsilon$, where $\varepsilon
\in (0,\beta-1)$; then $p \in (1,\beta)$ is equivalent to write
$1<\beta-\varepsilon<\beta$. Therefore
$\psi(p)=\psi(\beta-\varepsilon)=\varepsilon^{-\frac{\theta}{\beta-\varepsilon}}$
and the $G\psi$ norm of any measurable function $ f: Q \to \mathbb
R$ takes the well-known form
\begin{equation}\label{GLS-Lebesgue}
||f||_{G\psi} =
\sup_{0<\varepsilon<\beta-1}\varepsilon^{\frac{\theta}{\beta-\varepsilon}}||f||_{\beta-\varepsilon},
\end{equation}
introduced in \cite{Iwaniec-Sbordone} and denoted by
$||f||_{L^{\beta),\theta}} \, ( \, \beta>1, \, \theta\geq0$).

%
%\noindent  Let $\varepsilon \in (0,\beta-1)$ and replace \,$p$\,
%with $p-\varepsilon$ and $\beta$ with $p$ in the function $\psi(p)$.
%Therefore we have
%$\psi(p-\varepsilon)=\varepsilon^{-\frac{\theta}{p-\varepsilon}}$, \
%$\varepsilon\in (0,p-1)$, $p>1$, and the $G\psi$ norm of any
%measurable function $ f: Q \to \mathbb R$ takes the well-known form
%\begin{equation}\label{GLS-Lebesgue}
%||f||_{G\psi} =
%\sup_{0<\varepsilon<p-1}\varepsilon^{\frac{\theta}{p-\varepsilon}}||f||_{p-\varepsilon},
%\end{equation}
%introduced in \cite{Iwaniec-Sbordone} and denoted by
%$||f||_{L^{p),\theta}}$.

}
\end{remark}

\vspace{2mm}

 The Grand Lebesgue Spaces spaces have been studied (and applied) in huge numbers of works, see. e.g.
\cite{ForKozOstr_Lithuanian,Formica Ban Alg,Ahmed Fiorenza Formica
at all,anatriellofiojmaa2015,
anatrielloformicaricmat2016,BuldKoz1,caponeformicagiovanonlanal2013,Ermakov
etc. 1986,Fiorenza1,Fiorenza2,KozOstr1,Ostrovsky
01901,Ostrovsky1,Ostrovsky tail,Ostrovsky HIAT,Ostrovsky3}, etc.

We note that the $G\psi$ spaces are also interpolation spaces (the
so-called $\Sigma$-spaces), see \cite{Jawerth Milman 91}.
 \par

 The theory of GLS allows to investigate the exponential decreasing tail behavior of measurable functions, as well as
 estimates of its Orlicz norms, which is used in Probability theory, Functional Analysis, theory of Partial Differential
Equations, etc.

\vspace{4mm}

\section{Quasi-Grand Lebesgue Spaces. Some properties}\label{Section Quasi Grand}

\vspace{2mm}

Let again $(Q = \{ x\}, B,\mu) $ be a measure space with non-trivial
sigma-finite measure $ \mu$ \, and let
 $ \, \psi(p) = \psi_{a,b}(p), \ p \in (a,b)$, where now
\begin{equation} \label{a b}
 0 < a < b \le  1,
\end{equation}
be a continuous, strictly positive, numerical valued function in the
{\it open} interval $(a,b)$,  such that
$$
\inf_{p \in (a, b)} \psi(p) > 0.
$$

The Grand Lebesgue Space (GLS) $ \, G\psi = G\psi_{a,b}\,$ consists
of all the measurable functions $f: Q \to \mathbb R$ such that
\begin{equation} \label{quasi Grand Lebesgue}
||f||_{G\psi} = ||f||_{G\psi_{a,b}}  \stackrel{def}{=}  \sup_{p \in
(a,b)} \left\{ \ \frac{||f||_{p}}{\psi(p)} \ \right\}<\infty.
\end{equation}

 \noindent In contradistinction to the foregoing definition of ordinary GLS, given in Section \ref{Section_OrdinaryGLS}, these spaces are not Banach spaces. They are only quasi-Banach local
bounded $ \,F\,$ spaces, as proven in next Proposition \ref{Prop
quasi-norm}.

\vspace{0.2cm}

By definition, the following $\psi$\,-function
$$
\psi_f(p) \stackrel{def}{=} ||f||_p, \ \ \  p \in (a,b),
$$
is named {\it natural function} for the  real valued function $f$.

Evidently $f \in G\psi_f$ and
$$
||f||_{G\psi_f} = 1.
$$

\vspace{2mm}

\begin{example}{\rm \textbf{Natural function}.}

\vspace{2mm}

 {\rm Let $ f(x) = f_{\Delta, \delta, L}(x), \  x
\in (0,1)$, be the following function
$$
f(x) := x^{-\Delta} \ |\ln x|^{\delta} \ L(|\ln x|), \ \ \Delta =
{\const} > 1, \ \ \delta = {\const} \ge 0,
$$
 where $ L = L(y), \ y \in (0,\infty)$, is a positive continuous {\it slowly varying} function as $ y \to \infty $. \par

 The natural function for $f$ is
$$
\psi_f(p) = ||f||_p, \ \ \  p \in (0,1/\Delta).
$$

Obviously $f \in G\psi_{f}$ and
$$
||f||_{G\psi_{f}} = 1.
$$
 We observe, after some calculations, that the natural function
$ \psi_{f}(p) = ||f||_p \ $ for the function $f$ defined above has
the form, for $p \in (0, 1/\Delta)$,
$$
\psi_{f}(p) =||f||_p \, \asymp \, \Gamma^{1/p} (p \delta + 1) \  L
\left( \ \frac{1}{1 - p \, \Delta} \ \right) \ \cdot (1 - p
\Delta)^{-\Delta - 1/p},
$$
where $ \Gamma(\cdot) $ denotes the Euler's Gamma function.

It is easy to see that
$$
\lim_{p \to 0+}\psi_{f}(p)= \lim_{p \to 0+} ||f||_p = L(1) \exp
\left( \Delta - \delta\, \gamma \right),
$$
where $ \gamma $ is the Euler-Mascheroni constant
$$
\gamma = - \int_0^{\infty} \ln z \ e^{-z} \ dz \approx 0.57721566...
.
$$

}
\end{example}

\vspace{2mm}

 See also Example \ref{example tail} further.

 \vspace{2mm}

 \begin{remark}
{\rm More generally it is known that, for a function $ \ f: \ (0,1)
\to \mathbb R $,
 $$
 \lim_{p \to 0+} ||f||_p = \exp \left( \ \int_0^1 \ln |f(x)| \ dx \ \right).
 $$
 Namely, when $ \ f(x) > 0, \ x \in (0,1), \ $ then
$$
\int_0^1 f^p(x) dx = \int_0^1 \exp(p \ \ln f(x)) \ dx \ \buildrel{p
\to 0^+} \over {\sim} \  \left[ \ 1 + p \int_0^1 \ln f(x) \ dx \
\right],
$$
%\overset{p \to 0^+}{\sim}
so
$$
\lim_{p \to 0^+}||f||_p=%\left[ \ \int_0^1 f^p(x) dx \ \right]^{1/p}
\lim_{p \to 0^+}\left[ \ 1 + p \int_0^1 \ln f(x) \ dx \
\right]^{1/p} ß= \exp \left[   \ \int_0^1 \ln f(x) \ dx \ \right].
$$

}
 \end{remark}

 \vspace{2mm}

\begin{remark}\label{particular case Bandalayev}

{\rm  The particular case of these spaces, %namely the grand Lebesgue
%space $L^{p)}(0,1)$, \, $(0<p\leq 1)$,
was introduced in \cite{Bandaliyev}, where applications for the
study of Hardy operators were given.

 More precisely, let $Q=(0,1)$, $\displaystyle a=\frac{b}{2},0< b\leq 1$. Define $\psi(p)=(b-p)^{-\frac{1}{p}}, \ p \in \displaystyle\left(\frac{b}{2},b\right)$.
  Let $\varepsilon \in (0,b/2)$ and replace \,$p$\, with
$b-\varepsilon$ %and $b$ with $p$
in the function $\psi(p)$. Therefore we have
$\psi(p)=\psi(b-\varepsilon)=\varepsilon^{-\frac{1}{b-\varepsilon}}$,
\ $\varepsilon\in \left(0,\frac{b}{2}\right)$, $0<b\leq 1$, and the
$G\psi$ norm of any measurable function $ f: (0,1) \to \mathbb R$ is
given by
\begin{equation}\label{GLS-Lebesgue}
||f||_{G\psi} =
\sup_{0<\varepsilon<\frac{b}{2}}\varepsilon^{\frac{1}{b-\varepsilon}}||f||_{b-\varepsilon},
\end{equation}
so we recover the quasi-norm introduced in \cite{Bandaliyev},
denoted by $L^{b)}(0,1)$, \, $(0<b\leq 1)$.
 }

 \end{remark}

\vspace{2mm}

\begin{proposition}\label{Prop quasi-norm}

The Grand Lebesgue space $ G\psi=G\psi_{a,b}$,\, $0<a<b\leq 1$, is a
quasi-normed space, having quasi-triangle constant
\begin{equation} \label{trian a b}
C = C[\psi_{a,b}] = 2^{1/a - 1};
 \end{equation}
 i.e., let $f,g\in G\psi$, then
 $$
||f + g||_{G\psi} \le 2^{1/a - 1}  \, (\,||f||_{G\psi} +
||g||_{G\psi} \,).
$$

 \end{proposition}

\vspace{2mm}

\noindent \begin{proof}

\vspace{2mm}

 Let $ \ f,g \in G\psi$. The first two properties of the quasi-norm \eqref{axiom 1} and \eqref{axiom 3} trivially holds.  Now we have
$$
||f + g||_p^p = \int_Q |f(x) + g(x) |^p \ \mu(dx).
$$
 Applying inequality (\ref{ineq A B}) we have
$$
||f + g||_p^p \le \int_Q |f(x)|^p \ \mu(dx)  +  \int_Q |g(x) |^p \
\mu(dx) = ||f||_p^p +  ||g||_p^p \, ,
$$
so we get inequality (\ref{p norm}). Furthermore, using once again
(\ref{ineq A B}), we deduce
\begin{equation} \label{quasi est}
||f + g||_p \le 2^{1/p -1} \, (\, ||f||_p +  ||g||_p \,).
\end{equation}

Since $ \ f,g \in G\psi$ then $||f||_{G\psi} < \infty$ and
$||g||_{G\psi} < \infty$. Therefore
\begin{equation}\label{estimates norm f g}
||f||_p \le||f||_{G\psi}\cdot \psi(p), \ \ \ ||g||_p \le
||g||_{G\psi}\cdot \psi(p), \ \ p \in (a,b).
\end{equation}
 By \eqref{quasi est} and \eqref{estimates norm f g} it follows
\begin{equation*}
||f + g||_p \le 2^{1/p -1} (\, ||f||_{G\psi}+||g||_{G\psi}\,) \,
\psi(p), \ \ p \in (a,b).
\end{equation*}
Taking into account that $ \ p > a > 0 \ $, we get
$$
||f + g||_{G\psi} \le 2^{1/a - 1}  \, (\,||f||_{G\psi} +
||g||_{G\psi} \,),
$$
as desired.
\end{proof}

 The proof of completeness of the $G\psi_{a,b}$ spaces, \, $0<a<b\leq 1$, is given in the next
 Section. We also observe that $G\psi_{a,b}$ spaces are
 rearrangement invariant.

\vspace{2mm}

Let us discuss now the {\it fundamental function} for these spaces.
We suppose here that the measure $ \, \mu \, $ is not trivial
(non-zero) and atomless. \par Recall that the  fundamental function
of an arbitrary $ \, G\psi_{a,b} \, $ space has the form

\begin{equation} \label{fu fu}
\phi_{G\psi_{a,b}}(\delta) = \sup_{p \in (a,b)} \left\{ \
\frac{\delta^{1/p}}{\psi(p)} \ \right\}, \  \ \ \delta \in [0,
\mu(Q)].
\end{equation}

 \vspace{2mm}

This function plays a very important role in functional analysis,
theory of Fourier series, etc. It was investigated in detail
especially for ordinary GLS in \cite{Ostrovsky3}.

\vspace{2mm}

 {\it We extrapolate this notion} \eqref{fu fu} {\it one-to-one
also in the case of quasi-Grand Lebesgue Spaces, i.e. when } $0 < a
< b \leq 1$. \par

\vspace{2mm}
 Denote
$$
\overline{\psi} := \sup_{p \in (a,b)} \psi(p), \ \ \  \ \
\underline{\psi} := \inf_{p \in (a,b)} \psi(p).
$$
Recall that by definition $ \ \underline{\psi} > 0. \ $ \par

\vspace{4mm}

 We state the following simple  bilateral estimate for the fundamental function
for the values $ \ \delta \in [0,1]$.

\vspace{2mm}

\begin{proposition}\label{Prop fundamental function}
The fundamental function $\phi_{G\psi_{a,b}}(\delta)$ of the GLS
space $G\psi_{a,b}$,\, with $0<a<b\leq 1$, satisfies, for the values
$\delta \in [0,1]$,

\begin{equation} \label{estim fu fu}
\frac{\delta^{1/a}}{\overline{\psi}} \le \phi_{G\psi_{a,b}}(\delta)
\le \frac{\delta^{1/b}}{\underline{\psi}}.
\end{equation}
\end{proposition}

Of course, the left-hand side of  the estimate \eqref{estim fu fu}
has a sense only if $\overline{\psi} < \infty $.

 \vspace{4mm}

 \section{Completeness.}

\vspace{2mm}

 In this Section we prove that the quasi Grand Lebesgue spaces introduced in Section \ref{Section Quasi Grand} are complete.

 We will use a well known criterion of completeness of quasi-normed spaces, proved by Lech Maligranda, which we recall here for convenience of the reader.

%\newpage

\begin{theorem}{\rm \cite[Theorem 1.1]{Maligranda 1}}
\label{ThMaligranda}

A quasi-normed space $X=(X,||\cdot||)$ with a quasi-triangle
constant $C\geq 1$ is complete (quasi-Banach space) if and only if,
for every series such that
$$\sum_{k=1}^\infty C^k||x_k||<\infty,$$
we have $\displaystyle\sum_{k=1}^\infty x_k\in X $and
$$\left|\left|\sum_{k=1}^\infty x_k \right|\right|\leq C \sum_{k=1}^\infty C^k||x_k||.$$

\end{theorem}

\vspace{4mm}

Now we state our result.

\begin{theorem}\label{theorem completeness}
Let $(Q = \{ z\}, B,\mu) $ be a measure space with non-trivial
sigma-finite measure $ \mu$. The quasi-Grand Lebesgue Space
$\,G\psi_{a,b}, \ 0 < a < b \leq 1 $, builded on $(Q, B,\mu) $, is
complete with respect to the quasi-distance $ \ \rho(x,y) := ||x -
y||_{G\psi_{a,b}}$.

\end{theorem}

\vspace{4mm}

\noindent \begin{proof}

\vspace{2mm}

Let us consider first the space $ \ L_p, \ 0 < p < 1; \ $ the case $
\ p = 1 \ $ is trivial. We use Theorem \ref{ThMaligranda}, which we
transform as follows.

Note that for $L_p$ spaces, with $0<p<1$, the metric used is
$\rho(x,y)=||x-y||_p^p$, i.e. without the $p$-th root of the
integral, since in this case we have, for $x,y,w\in L_p$,
\begin{equation*}
\begin{split}
\rho(x,y)&=||x-y||_p^p=||(x-w)+(w-y)||_p^p\\
& \leq \int_Q|x-w|^p \,d\mu(z)
+\int_Q|w-y|^p\,d\mu(z)\\
& =\rho(x,w)+\rho(w,y)
\end{split}
\end{equation*}
with constant $C=1$.

\vspace{3mm}

Assume that for an arbitrary sequence $ \ \{x_k,\}, \ k =
1,2,\ldots, \ x_k \in L_p \ $ such that
\begin{equation} \label{p criterion}
S:= \sum_{k=1}^{\infty} ||x_k||^p_p < \infty,
\end{equation}
the sum
\begin{equation} \label{sum belongs}
x := \sum_{k=1}^{\infty} x_k
\end{equation}
belongs to the space $ \, L_p \, $ and moreover
\begin{equation} \label{sum belongs}
||x||_p^p \le \sum_{k=1}^{\infty} ||x_k||^p_p.
\end{equation}
Then actually the quasi-metric space $ \ (L_p,\rho) \ $ is complete.

\vspace{2mm}

So, let $ S < \infty$ or equally
$$
\sum_{k=1}^{\infty} \int_Q|x_k(z)|^p \mu(dz) < \infty.
$$
It follows immediately from the theorem of Beppo-Levi that the sum $
\ \sum_{k = 1}^{\infty} |x_k(z)|^p \ $ is integrable. Therefore
$$
||x||_p^p \le \sum_{k=1}^{\infty} ||x_k||^p_p = S < \infty.
$$

\vspace{2mm}

Let us consider now the case of the space $ \ G\psi = G\psi_{a,b}, \
0<a<b<1$. Let $ C = 2^{1/a - 1} \in (1, \infty)$ the quasi-triangle
constant (as seen in Proposition \ref{Prop quasi-norm}) and let $\{
x_k \}, \ k = 1,2,\ldots $, be any sequence of the elements of the
space $ G\psi_{a,b}$ such that
$$
\Theta := \sum_{k=1}^{\infty} C^k  \, ||x_k||_{G\psi_{a,b}}  <
\infty.
$$

We have

$$
\Theta = \sum_{k=1}^{\infty}  C^k  \sup_{p\in(a,b)} \left[ \
\frac{||x_k||_p}{\psi(p)} \right] < \infty,
$$

$$
\sum_{k=1}^{\infty} C^k ||x_k||_p  \le \Theta \ \psi(p).
$$

  As we yet know, the value $ \ x \ $ there exists as element of the space $ \ L_p(Q) \ $  for all the values $ \ p \in (a,b) \ $ and
moreover

$$
||x||_p \le C \Theta \, \psi(p).
$$
Consequently, $ x \in G\psi  $ and $ ||x||_{G\psi} \le C \Theta$.
\par

This completes the proof.

 \end{proof}

\vspace{4mm}

 \section{Dual space}\label{Section dual space}

\vspace{2mm}

 We suppose in this Section that the measure $ \, \mu \, $ is atomless and that the function $ \psi(p) = \psi_{a,b}(p), \ p\in (a,b)$,  \ $0 < a < b < 1$, is described as before:
 it is continuous and strictly positive in $(a,b)$ and  $\displaystyle \inf_{p \in
(a,b)} \psi(p) > 0$.\par

 We intend here to investigate the dual (conjugate) space of the quasi-GLS space $G\psi=G\psi_{a,b}$; more precisely, we will prove
 its essential absence. \par

 It is known that the dual space $\,L^*_p\, $ of the $\,L_p\,$ space, $p \in (0,1)$, is {\it degenerate:} it consists only of the null element (see \cite[pp. 36-37]{Rudin}). It is natural to expect that the same conclusion is also true for the quasi-Grand Lebesgue Spaces.

 \vspace{2mm}

\begin{theorem}\label{theorem dual}

The dual space $ \left[ \, G\psi_{a,b} \, \right]^*$ of the
quasi-Grand Lebesgue Space $\,G\psi_{a,b}, \ 0 < a < b < 1 $,
builded on a measure space  $(Q, B,\mu) $,  with atomless
sigma-finite measure $\mu$, is degenerate, that is
\begin{equation} \label{degenerateness}
 \left[ \, G\psi_{a,b} \, \right]^* = \{0\}.
\end{equation}

 \noindent As a consequence, the {\it associate} space  $\left[ \, G\psi_{a,b} \, \right]'$ is also degenerate.

 \end{theorem}

 \vspace{2mm}

 \noindent \begin{proof}

  We will follow mainly arguments like the ones described in the book of W. Rudin \cite[\S1.47, pp. 36-37]{Rudin}.
 Denote, for brevity, by $ \phi(\delta), \ \delta \in (0,1)$,
 the fundamental function $\phi_{G\psi_{a,b}}(\delta)$ for the space $ G\psi = G\psi_{a,b}$, defined in \eqref{fu fu}. \par

We claim that $G\psi$ contains no convex open sets, other than
$\emptyset$ and $G\psi$.

To prove this, suppose $ V \ne \emptyset$ is an arbitrary open and
convex set in $ G\psi$. Assume, without loss of generality, that
$0\in V $. Then, for
 some positive finite number $r$, \ $B_r \subset V, \ $ where $ \ B_r \ $ is the centered open ball in $
 G\psi$ \,:
$$
B_r = \{g \in G\psi \, : \, ||g||_{G\psi} < r \ \}.
$$
 Let us pick an arbitrary non-zero function $ f\in G\psi$ such that $||f|| = ||f||_{G\psi} \in
 (0,\infty)$.

 \vspace{2mm}

 \noindent Let $n\in \mathbb N, \, n\geq 2$, and introduce the following  measurable {\it partition} $ \ P = P_n \ $ of the  whole set $Q
 $ \,:
\begin{equation} \label{partition cap}
P_n  =  \{ \ A_n(i) \ \},  \ \ i = 1,2,\ldots,n; \  \ i \ne j \ \ :
\ \ A_n(i) \cap A_n(j)  = \emptyset,
\end{equation}
\begin{equation} \label{partition cup}
\bigcup_{i=1}^n A_n(i) = Q,
\end{equation}
 and a corresponding double sequence of (measurable) functions

$$
g_{i,n} = n \cdot f \cdot I_{A_n(i)},
$$
 where $I_{A_n(i)} $ denotes the indicator function of the sets $A_n(i)$,  such that

$$
||g_{i,n}|| =  n \cdot \phi(1/n) \cdot||f||.
$$

\vspace{1mm}

 \noindent By the right-hand side of (\ref{estim fu fu})
in Proposition \ref{Prop fundamental function} we have
\begin{equation} \label{ upp estim fu fu}
%\phi_{G\psi_{a,b}}(1/n)
\phi(1/n)\le \frac{n^{-1/b}}{\underline{\psi}},
\end{equation}
 and we get
$$
\forall i = 1,2,\ldots,n \ \ \Rightarrow \ \ ||g_{i,n}|| \le
\frac{n^{1 - 1/b}}{\underline{\psi}} \ \cdot ||f|| \le r, \  \ n >2,
$$
 since $ \ n^{1- 1/b} \to 0$ as $n \to \infty $.  On the other words, for all the large values $n$, \ $n>2$,

$$
\forall i = 1,2, \ldots n \ \  \Rightarrow \ \ g_{i,n} \in B_r
\subset V.
$$
On the other hand,

$$
f = \frac{1}{n}\sum_{i=1}^n g_{i,n}.
$$
As long as the set $ \ V \ $ is convex, we deduce from the last
equality that $ \ f \in V. \ $ \par

\vspace{2mm}

Thus, an arbitrary non empty convex set in the space $ \ G\psi \ $
coincides with the whole space $ \ G\psi\ $. The lack of convex open
sets implies, as consequence, that the dual space of $G\psi$ is
degenerate, that is $[G\psi]^*  = \{0\} $ (see \cite[pp.
36-37]{Rudin}).

 \end{proof}

\vspace{2mm}

 Theorem \ref{theorem dual}  may be generalized as follows.

\begin{theorem}
 Let $ Y = (Y,||\cdot||_Y) $ be an arbitrary quasi-Banach rearrangement
 invariant function space builded over an atomless sigma-finite measure space $ (Q,B,\mu)
 $. Let $ f\in Y$ be any non-zero function such that $0 < ||f||_Y <
\infty$ and define the following generalization of the notion of the
fundamental function (weighted version):
$$
\phi_{Y,f}(\delta) \stackrel{def}{=} \sup_{A \in B, \ \mu(A) \le
\delta} ||f \cdot I_A||_Y, \ \ \ \delta \in (0,\infty).
$$
Evidently,
$$
\phi_{Y,f}(\delta) \le \esssup_{v \in Q} |f(v)| \cdot
\phi_{Y}(\delta),
$$
where $\phi_Y$ is the ordinary fundamental function.

Consider some partition of the whole space $ Q $ containing $ \, n
\, $ elements

$$
P_n = P_n(\, \{A_n(i)\} \, ), \ \ \ i = 1,2,\ldots,n; \ \ \ n =
2,3,\ldots,
$$
and denote the collection of such partitions with $S(n) := \{ P_n
\}$. Define the following variables:
\begin{equation} \label{def H}
H[P_n,f](Y) \stackrel{def}{=} \inf_{ \{A_n(i)\} \in P_n } \max_{i =
1,2,\ldots,n} \phi_{Y,f}(\mu(A_n(i)),
\end{equation}

\begin{equation} \label{def overline H}
\overline{H}_n [Y](f) \stackrel{def}{=} \inf_{P_n  \in S(n)}
H[P_n,f](Y).
\end{equation}

 We assert that, if for an arbitrary non-zero function $ f\in Y $, we
 have

\begin{equation} \label{condition H}
\lim_{n \to \infty}  n \ \overline{H}_n[Y](f) = 0,
\end{equation}
then the dual space $Y^* $ of \, $ Y$ is degenerate, i.e. $ Y^* =
\{0\} $.

\end{theorem}

\noindent\begin{proof}

Suppose $ V\subset Y$,\ $ V \ne Y$, is an arbitrary centered (i.e.
$0\in V$) non-trivial convex set and suppose that $ \exists\, r \in
(0,\infty)$ such that $B_r  \subset V$, where
$$
B_r = \{g \in Y \, : \, ||g||_{Y} < r \ \}.
$$
Let $ f\in Y$ be any non-zero function such that $0 < ||f||_Y <
\infty$.

%Define the following generalization of the notion of the fundamental
%function, weighted version:
%$$
%\phi_{Y,f}(\delta) \stackrel{def}{=} \sup_{A \in B, \ \mu(A) \le
%\delta} ||f \cdot I_A||_Y, \ \ \ \delta \in (0,\infty).
%$$
%
% Evidently,
%$$
%\phi_{Y,f}(\delta) \le \esssup_{v \in Q} |f(v)| \cdot
%\phi_{Y}(\delta),
%$$
%where $\phi_Y$ is the ordinary fundamental function.
%
%\vspace{2mm}
%
% We consider some partition of the whole space $ Q $ containing $ \, n \, $ elements
%
%$$
%P_n = P_n(\, \{A(i,n)\} \, ), \ \ \ i = 1,2,\ldots,n; \ \ \ n =
%2,3,\ldots,
%$$
%and denote the collection of such partitions as $S(n) := \{ P_n \}$.
%Define the following variables:
%
%\begin{equation} \label{def H}
%H[P_n,f](Y) \stackrel{def}{=} \inf_{ \{A(i,n)\} \in P_n } \max_{i = 1,2,\ldots,n} \phi_{Y,f}(\mu(A(i,n)),
%\end{equation}
%
%\begin{equation} \label{def overline H}
%\overline{H}_n [Y](f) \stackrel{def}{=} \inf_{P_n  \in S(n)} H[P_n,f](Y).
%\end{equation}
%
% We assert that, for an arbitrary non-zero function $ f\in Y $, we
% have
%
%\begin{equation} \label{condition H}
%\lim_{n \to \infty}  n \ \overline{H}_n[Y](f) = 0.
%\end{equation}

 As in Theorem \ref{theorem dual}, let $ \  n \ge 2 \  \ $ be a certain "great" natural number and let
 $ \ P_n = \{ A_n(i) \}  \ $  be an optimal partition in the sense of \eqref{def H}, \ \eqref{def overline H}:
$$
H[P_n,f](Y) = \overline{H}_n [Y](f).
$$
We pick
$$
g_{i,n} = n \cdot f \cdot I_{A_n(i)},
$$
 so that, as $ \ n \to \infty \ $,
$$
||g_{i,n}||_Y  = n \ \phi_{Y,f} \ (\mu(A_n(i)) \le n \
\overline{H}_n [Y](f) \to 0,
$$
by virtue  of the condition \eqref{condition H}.

\noindent Therefore, for $n$ sufficiently \lq\lq large\rq\rq,
$g_{i,n}\in B_r\subset V, \ \forall\, i=1,\ldots,n$.
\par

 We observe that
$$
f = \frac{1}{n}\, \sum_{i=1}^n g_{i,n},
$$
i.e. the (arbitrary) function $ \ f \ $ belongs to the convex linear shell of the elements $ \ g_{i,n} \ $ from the set $ \ V. \ $
 As long as  the set $ \ V \ $ is convex, we conclude $ \ f \in V, \ $ therefore $ \ V = Y, \ $ in contradiction.

 Thus, an arbitrary non empty convex set in the space $ Y $
coincides with the whole space $ Y$. For the same arguments at the
end of the previous Theorem, we conclude that the dual space of $Y$
is degenerate, that is $Y^* = \{0\} $.

\end{proof}

\begin{remark}
{\rm As a particular case of $G\psi$ spaces, thanks to Remark
\ref{particular case Bandalayev}, we have that the dual of the grand
Lebesgue space $L^{b)}(0,1), \ 0<b< 1$, described in
\cite{Bandaliyev}, is $[L^{b)}(0,1)]^*=\{0\}$.}

\end{remark}

For the interested reader, the dual of the grand Lebesgue space
$L^{b)}(0,1)$, \ $b>1$, has been described in \cite{Fiorenza1} (see
also
\cite{caponefiorenzajfsa2005,CaponeFormica_dual,Anatriello_Formica_Giova_JMAA2017}),
while for the dual of the general (ordinary) GLS spaces
$G\psi=G\psi_{\alpha,\beta}$, \ $1\leq \alpha<\beta\leq\infty$, see
\cite{Ostrovsky_bilateral_small2009,Ostrovsky_dual2017}.

  \vspace{4mm}

 \section{Other  results.}

\vspace{2mm}

 %\begin{center}
%
% {\sc  Connections between tail behavior and quasi Grand Lebesgue Norm.}
%
% \end{center}

 \subsection{Connections between tail behavior and quasi Grand Lebesgue Norm}

Let $ \ (Q, B,\mu) \ $ be a non-trivial measure space with
sigma-finite measure $\mu$. We define, for an arbitrary measurable
function $ \ f: Q \to \mathbb R \ $, its so-called tail function

 $$
 T_f(u) \stackrel{def}{=} \mu \{ \ z: \ |f(z)| \ge u \  \}, \  \ u > 0.
 $$

It follows from Tchebychev inequality, analogously as done in
\cite{Ostrovsky HIAT},

$$
T_f(u) \le \frac{||f||_p^p}{u^p} \le \frac{[||f||_{G\psi}]^p \
\psi^p(p)}{u^p},
$$
therefore

\begin{equation} \label{norm tail}
T_f(u) \le \inf_{p \in (a,b)} \left\{ \ \frac{||f||_p^p}{u^p}
\right\}  \le  \inf_{p \in (a,b)} \ \left\{ \frac{[||f||_{G\psi}]^p
\ \psi^p(p)}{u^p} \ \right\},
\end{equation}
where $0<a<b\leq 1$.

\vspace{2mm}

Conversely, introduce the following $ \, \psi \, $ function

$$
\psi_0(p) := \left[ p \int_0^{\infty} \ u^{p-1} T_f(u) \ du \
\right]^{1/p} = ||f||_p ,
$$
If the last integral converges in some interval $ \ (a,b), \
0<a<b\leq1 \ $, then
\begin{equation} \label{tail norm}
f \in G\psi_0, \ \ \  ||f||_{G\psi_0} = 1.
\end{equation}

\vspace{4mm}

 Let us show an example.

\vspace{2mm}

\begin{example}\label{example tail}

{\rm Define the following tail function

$$
T^{b,\gamma, L}(x) = x^{-b} \ (\ln x)^{\gamma} \ L(\ln x), \ \ \ x
\ge e,
$$
where  $ \ b \in (0,1), \ \gamma = {\const} \ge 0, \ \mu(Q) = 1, \ L
= L(y), \  y\geq 1$, is a positive continuous {\it slowly varying}
function as $ y \to \infty $. We deduce after some calculations,
alike ones in  \cite{Ostrovsky 01901}, that if a measurable function
(random variable) $\ \xi \ge 1 \ $, defined on some probability
space, is such that
%&&
\begin{equation}\label{estimate tail}
T_{\xi}(x) \le C \ T^{b,\gamma, L}(x), \ \ \ x \ge e,
\end{equation}
then, for $  p \in (0,b), \ b \in (0,1), \ x\geq e$
\begin{eqnarray*}  %\label{tail p-beta}
\begin{split}
||\xi||_p=& \left[ p \int_0^{\infty} \ x^{p-1} T_{\xi}(x) \ dx \
\right]^{1/p} \\
\leq \, & C  \left[ p \int_e^{\infty} \ x^{p-1} x^{-b} \ (\ln
x)^{\gamma}
\ L(\ln x) \ dx \ \right]^{1/p} \\
\leq \, & C  \left[ p \int_1^{\infty} \ x^{p-1} x^{-b} \ (\ln
x)^{\gamma}
\ L(\ln x) \ dx \ \right]^{1/p} \\
= \,& C\left[ p \int_0^{\infty} \ e^{-y(b-p)}\, y^\gamma L(y)\, dy  \right]^{1/p}\\
= \,& C \left[ p \int_0^{\infty} \ e^{-z}\, (b-p)^{-\gamma}\,
z^\gamma L\left(\frac{z}{b-p}\right) (b-p)^{-1}\, dz  \right]^{1/p}  \\ \\
\sim \,&  C \, \left[ p (b-p)^{-1} (\beta-p)^{-\gamma}
L\left(\frac{1}{b-p}\right)\int_0^\infty e^{-z} z^\gamma\, dz\right]^{1/p} \hspace{1cm} ({\rm as} \ p\to b^-) \\
\\
=\, & C \left [p \,(b-p)^{-(\gamma+1)}L\left(\frac{1}{b-p}\right)
\Gamma(\gamma+1)\right]^{1/p}\\ \\
 = \, & C_1(b,\gamma,L) \ (b - p)^{-(\gamma
+ 1)/p}
\ L^{1/p}(1/(b - p)) \\ \\
\asymp \ & C_2(b,\gamma,L) \ (b - p)^{-(\gamma + 1)/b} \
L^{1/b}(1/(b - p))
\end{split}
\end{eqnarray*}
%
%\begin{equation} \label{tail beta}
% C_2(\beta,\gamma,L) \ (\beta - p)^{-(\gamma + 1)/\beta} \ L^{1/\beta}(1/(\beta - p)) =: C_2 \psi_{\beta,\gamma, L}(p).
%\end{equation}

In conclusion
\begin{equation}\label{tail p-beta}
||\xi||_p\leq C_2(b,\gamma,L) \ (b - p)^{-(\gamma + 1)/b} \
L^{1/b}(1/(b - p))\ =: \ C_2\,\psi(p)
\end{equation}

We state that the random variable $ \, \xi \, $ belongs to the
quasi-Grand Lebesgue Space $ \, G\psi\,$, where
$$
\psi(p)=(b - p)^{-(\gamma + 1)/b} \ L^{1/b}(1/(b - p)), \ \ \
p\in(0,b), \ \ \ b \in (0,1).
$$

\vspace{3mm}

  Inversely, suppose that the tail estimate
$$
||\xi||_p\leq C_1(b,\gamma,L) \ (\beta - p)^{-(\gamma + 1)/p} \
L^{1/p}(1/(b - p))
$$
holds for some random variable $ \xi$, under the restrictions $ p
\in (0,b), \ b \in (0,1)$.
 We obtain, by Tchebychev's inequality,
\begin{equation*}
\begin{split}
T_\xi(x) \le \frac{||\xi||_p^p}{x^p}\leq \  C_1^p \,x^{-p}\ (b-
p)^{-(\gamma + 1)} \ L(1/(b - p)),
\end{split}
\end{equation*}
and, choosing \ $ \displaystyle p := b - \frac{c}{\ln x}, \  x \ge e
$, we get the following tail estimate
\begin{equation} \label{Gpsi estim}
T_{\xi}(x) \le \  C_3 \, x^{-p} \ (\ln x)^{\gamma + 1} \ L\left(\ln
x\right), \ \ \ \ \  x \ge e.
\end{equation}

\vspace{3mm}

 It is very important to notice that there is a "gap" $ \ (\ln x)^1 \ $ between the estimates \eqref{Gpsi estim} and \eqref{estimate tail}.
This "gap" is essential, see correspondent examples in
\cite{Ostrovsky HIAT,Ostrovsky 01901}. \par
 This example remains true still for the values $ \ p \in (0,1) $.
}

\end{example}

 \vspace{4mm}

 \subsection{Boyd indices of the quasi-Grand Lebesgue spaces $\bm{G\psi_{a,b}}$}

 %{\sc Boyd's  indices of our space.} \par

%\vspace{4mm}

  Let now $  \ Q = \mathbb R_+ = ( 0,\infty)$ and $d\mu(z) = dz$. We provide the values of the Boyd indices for the quasi-GLS space $ \ G\psi=G\psi_{a,b} \ $ on $(0,\infty)$,
  having parameters $a,b$ such that $0 < a < b \le 1$.

Boyd indices play an important role in the theory of interpolation
of operators and in Fourier Analysis.

We recall the definition (see, e.g.,
 \cite{Bennet}). Denote
 $\sigma_s:G\psi_{a,b}\to G\psi_{a,b}$ the dilatation operators
 given by
 $$\sigma_sf(x)=f(x/s), \ \ \ s>0.$$
 The Boyd indices $\gamma_1$ and $\gamma_2$ of the quasi- GLS space $ G\psi$, with $\gamma_1\leq\gamma_2$, are defined by
\begin{equation}
\gamma_1=\gamma_1[G\psi]= \lim_{s\to
0^+}\frac{\log||\sigma_s||_{G\psi\to G\psi}}{\log s}, \ \ \ \ \
\gamma_2=\gamma_1[G\psi]= \lim_{s\to
\infty}\frac{\log||\sigma_s||_{G\psi\to G\psi}}{\log s}
 \end{equation}

 Following the same arguments concerning the classical GLS spaces (see
 \cite{liflyandostrovskysirotaturkish2010}), we state that
$$
\gamma_1 = \gamma_1[G\psi] = \frac{1}{b}, \ \ \ \  \  \gamma_2 =
\gamma_2[G\psi] = \frac{1}{a}.
$$

\vspace{3mm}

 We recall that the Boyd indices for the particular case of the
classical GLS spaces, namely for the space $L^{b)}(0,1), \ b>1$,
were computed in \cite{Formica_Giova_Boyd_indices2015}.

 \vspace{4mm}

%{\sc Contraction principle in quasi metric space.}
%
%\vspace{4mm}

\subsection{Contraction principle in quasi metric spaces}

 Let $ \ (X,d) \ $ be a {\it complete, closed}, non-trivial quasi-metric space, for instance, a quasi-Banach space, equipped with the
 quasi-distance function $ \ d = d(x,y), \ x,y \in X $,  such that
\begin{equation} \label{quasi dist}
\exists K = {\const} \in (1,\infty)   \ : \ d(x,y) \le K [d(x,z)+
d(z,y)], \ \ x,y,z \in X.
 \end{equation}
 The map $ \ f: X \to X \ $ is said to be a {\it contraction} with
parameter  $ \ \alpha = \const \in (0,1), \ $ iff the Lipschitz
condition is satisfied:
\begin{equation} \label{contract}
d(f(x), f(y) ) \le \alpha \ d(x,y).
\end{equation}

\vspace{4mm}

\begin{theorem}\label{th contraction}

Let $ (X,d)$ be a {\it complete, closed}, non-trivial quasi-metric
space and $ \ f: X \to X \ $ a contraction, with constant $\alpha$
in \eqref{contract} satisfying
 \begin{equation} \label{cond alpha}
 \alpha < \frac{1}{K^2},
 \end{equation}
 where $K$ is the constant in \eqref{quasi dist}.

 Then $f$ admits a unique fixed point $x^*\in X$, i.e. $x^*$ is the
 unique solution of the equation
\begin{equation} \label{fixed point eq}
f(x) = x,
\end{equation}
and $x^*$ may be obtained as limit, as $ n \to \infty $, of the
iterations (recursion)
\begin{equation} \label{recursion}
x_{n+1} = f(x_n), \ \ \ n = 0,1,2,\ldots,
\end{equation}
where $x_0$ is an arbitrary point in $X $. Herewith

\begin{equation} \label{error estim}
d(x^*,x_n) \le \frac{K \,\alpha^n}{1-\alpha \, K^2} \ d(x_0,x_1).
\end{equation}

\end{theorem}

\vspace{4mm}

\noindent \begin{proof}

The proof is a simple generalization of the one presented by R. S.
Palais in \cite{Palais}. First of all note that, for $x,y,z,v\in X$,
from \eqref{quasi dist} we have
\begin{equation} \label{K K2 K2}
d(x,y) \le K d(x,z) + K^2 d(z,v) + K^2 d(v,y),
\end{equation}
therefore, for arbitrary elements $ y_1, \ y_2 \in X $, by
\eqref{contract} follows
\begin{eqnarray*}
\begin{split}
d(y_1,y_2) & \le K d(y_1, f(y_1)) + K^2 d(f(y_1), f(y_2) ) + K^2
d(y_2, f(y_2)) \\ \\
& \leq  K d(y_1, f(y_1)) + \alpha K^2 d(y_1,y_2) + K^2 d(y_2,
f(y_2)).
 \end{split}
\end{eqnarray*}
Consequently,
\begin{equation} \label{key d estim}
d(y_1,y_2) \le \frac{K d(y_1,f(y_1)) + K^2 d(y_2, f(y_2))}{1 -
\alpha K^2}\, ,
\end{equation}
since $ \alpha K^2 < 1$.

\vspace{2mm}

Choosing in (\ref{key d estim})
$$
y_1 := x_n, \ \  \ y_2 := x_m, \ \ \ m > n \ge 1
$$
and taking into account the inequalities
$$
d(x_n,x_{n+1}) \le d_0 \alpha^n, \ \ \  d(x_m,x_{m+1}) \le d_0
\alpha^m,
$$
where $ \ d_0 = d(x_1,x_0) $, we get
\begin{eqnarray}\label{fund 1}
\begin{split}
d(x_n,x_m) & \le \frac{K d(x_n,x_{n+1}) + K^2 d(x_m, x_{m+1})}{1 -
\alpha K^2}\\ \\
&\leq \frac{ K d_0 \alpha^n + K^2 d_0 \alpha^m}{1 - \alpha K^2}, \ \
\ 1 \le n < m.
\end{split}
\end{eqnarray}

%\begin{equation} \label{fund 2}
%\frac{ K d_0 \alpha^n + K^2 d_0 \alpha^m}{1 - \alpha K^2}, \ 1 \le n < m.
%\end{equation}

The last estimate in (\ref{fund 1}) denotes that the the sequence $
\{x_n\}$ is fundamental (Cauchy sequence). Since the space $(X,d)$
is complete, there exists the limit
$$
x^* = \lim_{n \to \infty} x_n \in X.
$$

\noindent Finally, passing to the limit as $ \ m \to  \infty \ $ in
\eqref{fund
 1}, since $\alpha\in(0,1)$, we conclude
$$
d(x_n, x^*) \le \frac{K \, d_0 \, \alpha^n}{1-\alpha K^2}\,.
$$
This completes the proof.

\end{proof}

\vspace{3mm}

\begin{remark}

 {\rm  As we know, the condition (\ref{quasi dist}) is satisfied for the quasi Grand Lebesgue Spaces
 $ \ G\psi_{a,b}, \ 0 < a < b \le 1 \ $ with respect to the quasi-distance function $ \  d(f,g) = ||f - g||_{G\psi_{a,b}}$,
  with constant $ \ K =   2^{1/a - 1} $. }
\end{remark}

\vspace{3mm}

 Assume now that the quadrilateral inequality
\begin{equation} \label{K K K}
\exists\, K\in[1,\infty) \ : \ d(x,y) \le K [d(x,z) + d(z,v) +
d(v,y)],
\end{equation}
holds instead of the one in \eqref{K K2 K2}. Then Theorem \ref{th
contraction} may be reformulated as follows.

\vspace{4mm}

\begin{theorem}\label{th contraction KKK}

Let $ (X,d)$ be a {\it complete, closed}, non-trivial quasi-metric
space and $ \ f: X \to X \ $ a contraction, with constant $\alpha$
in \eqref{contract} such that
 \begin{equation} \label{cond alpha K}
 \alpha < \frac{1}{K},
 \end{equation}
 where $K$ is the constant in \eqref{K K K}.

 Then $f$ admits a unique fixed point $x^*\in X$, i.e. $x^*$ is the
 unique solution of the equation
\begin{equation} \label{fixed point eq}
f(x) = x,
\end{equation}
and $x^*$ may be obtained as limit, as $ n \to \infty $, of the
iterations (recursion)
\begin{equation} \label{recursion}
x_{n+1} = f(x_n), \ \ \ n = 0,1,2,\ldots,
\end{equation}
where $x_0$ is an arbitrary point in $X $. Herewith

\begin{equation} \label{error estim K}
d(x^*,x_n) \le \frac{K \,\alpha^n}{1-\alpha \, K} \ d(x_0,x_1).
\end{equation}

\end{theorem}

\vspace{3mm}

\begin{remark}
{\rm  Note that the quadrilateral inequality \eqref{K K K} is
satisfied for the quasi-Banach (Lebesgue-Riesz)
 spaces $ \ L_p, \ p \in (0,1)$, with constant $ \ K = 3^{1/p - 1}  \ $ as well as for the quasi-GLS spaces $ \ G\psi_{a,b}, \ 0 < a < b \le 1; \ $
in this case
$$
K=K[G\psi_{a,b}] := 3^{1/a - 1}.
$$
The last constant $ \, 3^{1/a - 1} \, $ is lesser than its old value
$ \ K^2 = \left[2^{1/a - 1}  \right]^2 $.}
\end{remark}

\vspace{3mm}

Some results concerning fixed point theorem in quasi-Banach spaces,
with applications to integral equations, was given in
\cite{Hussain}.

 \vspace{4mm}

 \section{Estimates for operators in quasi-Grand Lebesgue Spaces}

\vspace{2mm}

Let $\, U: L_p \to L_p, \  p \in (a,b)$, \, $ 0<a<b\leq1$, be a
bounded operator, not necessarily linear or sublinear, acting from
the quasi-Lebesgue-Riesz space $ \, L_p \, $ into itself: $
\,\forall \,f \in L_p \ \exists\, U[f] \in L_p \ $ and
\begin{equation} \label{Lp Lp}
 ||U[f]||_p \le \Theta(p) \ ||f||_p, \ \ \ p \in (a,b),
\end{equation}

$$
\forall p \in (a,b) \ \Rightarrow \ \Theta(p) < \infty.
$$

We can and will understood as the value of $ \ \Theta(p) \ $ its
minimal value, indeed

$$
\Theta(p) \stackrel{def}{=} \sup_{0 \ne f \in L_p} \left\{ \  \frac{||U[f]||_p}{||f||_p} \ \right\}.
$$

 Let also $ \ G\psi = G\psi_{a,b}, \ 0 < a < b \le 1 \ $  be again the quasi-Grand Lebesgue Space and let $ \ g \in G\psi_{a,b}, \ $
then

$$
||g||_p \le ||g||_{G\psi} \cdot \psi(p), \ \ \ p \in (a,b).
$$

 We conclude, using \eqref{Lp Lp},

\begin{equation} \label{Gpsi Gpsi}
||U[g]||_p \le \Theta(p) \ \psi(p) \ ||g||_{G\psi}.
\end{equation}

 On the other words, if we introduce a new $ \ \psi \ $ function $ \ \Psi[\Theta] (p) := \Theta(p) \
 \psi(p)$, we have
\begin{equation} \label{psi Theta}
||U[g]||_{G\Psi[\Theta]} \le ||g||_{G\psi}.
\end{equation}

\vspace{4mm}

\begin{proposition}

 \vspace{3mm}

\begin{equation} \label{Gpsi GTheta}
||U[\cdot]||\,_{\{G\psi  \to G\Psi[\Theta] \, \}}  = 1.
\end{equation}

\end{proposition}

\vspace{3mm}

\noindent \begin{proof}

The {\it upper estimate} in \eqref{Gpsi GTheta}, i.e. $ \
||U[\cdot]||\,_{\{G\psi  \to G\Psi[\Theta]\, \}} \le 1 $, follows
immediately from  \eqref{psi Theta} and  from the direct definition
of the norm in  the quasi Grand Lebesgue Spaces.\par

The {\it lower estimate} may be proved quite analogously  to the one
for the ordinary GLS, see \cite{Ostrovsky4}, \cite[Theorem
2.1.]{Ostrovsky5}.
\end{proof}

\vspace{2mm}

 Note, in addition, that the weighted estimates for Hardy and Hausdorff operators in quasi Banach spaces $L^{p)}$ and $L^p$ was obtained in
 \cite{Bandaliyev,Burenkov1,Burenkov2, Azzouz}.

\vspace{4mm}

 \section{Concluding remarks.}

\vspace{2mm}

\hspace{3mm} {\bf A.} Note that  the case $ \ K = 1 \ $  in Theorems
\ref{th contraction} and \ref{th contraction KKK} corresponds to the
classical contraction principle.
\par

\vspace{3mm}

  {\bf B.} \ Perhaps, the results obtained in this paper may be applied, for istance, in the investigation of the sums of random variables
having heavy tails of distributions, in the spirit of articles
\cite{OstrSir heavy,Ostrovsky tail}, etc., and, as a consequence, in
the Monte-Carlo method (\cite{Frolov,Grigorjeva}).

\vspace{3mm}

\vspace{0.5cm} \emph{Acknowledgement.} {\footnotesize The first
author has been partially supported by the Gruppo Nazionale per
l'Analisi Matematica, la Probabilit\`a e le loro Applicazioni
(GNAMPA) of the Istituto Nazionale di Alta Matematica (INdAM) and by
Universit\`a degli Studi di Napoli Parthenope through the project
\lq\lq sostegno alla Ricerca individuale\rq\rq .\par

\end{document}